\documentclass[11pt, english]{article}
\usepackage[top=1in, bottom=1in, left=1in, right=1in]{geometry} 
\usepackage[T1]{fontenc}
\usepackage{babel}

\usepackage[pdftex]{graphicx}
\usepackage{amssymb}
\usepackage{amsmath}
\usepackage{textcomp}
\usepackage{amsthm}
\usepackage{algorithm}
\usepackage{algorithmic}
\usepackage{subfigure}

\newtheorem{thm}{Theorem}

\newtheorem{lem}{Lemma}

\title{Small $H$-coloring problems for bounded degree digraphs}
\author{Pavol Hell \footnote{Partially supported by a NSERC discovery grant}\\ School of Computing Science,\\ Simon Fraser University,\footnote{The facilities of the IRMACS center at SFU are gratefully acknowledged, where most of this research was done} \\ Burnaby, B.C., Canada - V5A 1S6, \\ \texttt{pavol@cs.sfu.ca}
        \and Aurosish Mishra\footnote{Supported by the MITACS Globalink Internship Program for undergraduate students} \\ Department of Computer Science, \\ Cornell University, \\ Ithaca, NY, USA - 14853, \\ \texttt{aurosish@cs.cornell.edu}}
\date{}

\begin{document}
\maketitle

\begin{abstract}
An NP-complete coloring or homomorphism problem may become polynomial time solvable when restricted to graphs with degrees bounded by a small number, but remain NP-complete if the bound is higher. For instance, $3$-colorability of graphs with degrees bounded by $3$ can be decided by Brooks' theorem, while for graphs with degrees bounded by $4$, the $3$-colorability problem is NP-complete. We investigate an analogous phenomenon for digraphs, focusing on the three smallest digraphs $H$ with NP-complete $H$-colorability problems. It turns out that in all three cases the $H$-coloring problem is polynomial time solvable for digraphs with degree bounds $\Delta^{+} \leq 1$, $\Delta^{-} \leq 2$ (or $\Delta^{+} \leq 2$, $\Delta^{-} \leq 1$). On the other hand with degree bounds $\Delta^{+} \leq 2$, $\Delta^{-} \leq 2$, all three problems are again NP-complete. A conjecture proposed for graphs $H$ by Feder, Hell and Huang states that any variant of the $H$-coloring problem which is NP-complete without degree constraints is also NP-complete with degree constraints, provided the degree bounds are high enough. Our study is the first confirmation that the conjecture may also apply to digraphs.
\end{abstract}

\section{Introduction}

Graph coloring problems arise naturally in several contexts of both theoretical and applied nature. Be it scheduling events, solving pattern matching problems, or allocating registers to processes, many constraint satisfaction problems can be easily modeled in the graph coloring setting. A slightly more general problem is that of $H$-coloring.

Consider a fixed graph $H$. A {\em homomorphism} $f: G \rightarrow H$ is a mapping $f: V(G) \rightarrow V(H)$ such that $f(u)f(v)$ is an edge of $H$ for each edge $uv$ of $G$. An $H$-{\em coloring} of $G$ is a homomorphism $G \rightarrow H$. The definition is formally the same for digraphs $G, H$: since the edges are directed, the homomorphism $f$ is a mapping that preserves both the edges and their direction.

The $H$-{\em coloring problem} is the decision problem which asks the following question:
\begin{center}
 \textit{Given a graph $G$, is it possible to $H$-color $G$ $?$}
\end{center}

The complexity of the $H$-coloring problem for graphs $H$ has been widely studied \cite{hombook}. It is shown in \cite{HellNesetril1} that the $H$-coloring problem is polynomial time solvable if $H$ is bipartite or contains a loop, and is NP-complete otherwise. A number of variants of this basic family of problems have been considered, and in \cite{HellHuang} the authors have set up a framework for these variants. One of the parameters in this framework is a restriction to graphs with a given upper bound on the vertex degrees. For graphs with degrees bounded by $3$, some $H$-coloring problems that are NP-complete in general, become polynomial time solvable. However, these problems tend to be NP-complete again when the degree bound is $4$. Such is, for instance, the situation with $H = K_3$, i.e. with $3$-colorings. The theorem of Brooks ensures that a connected graph with degrees at most $3$, other than $K_4$, is $3$-colorable \cite{brooks}. Thus $3$-colorability of such graphs is decidable in polynomial time. However, it is known that the $3$-colorability of graphs with degrees at most $4$ is NP-complete \cite{Holyer}. (Our results also imply this fact, see Section 4.) Based on additional evidence of this kind, Feder, Hell and Huang conjectured that any variant of the $H$-coloring problem which is NP-complete without degree constraints is also NP-complete with degree constraints, provided the degree bound is high enough \cite{FederHellHuang}.

\begin{figure}
\centering
\includegraphics[scale = 0.6]{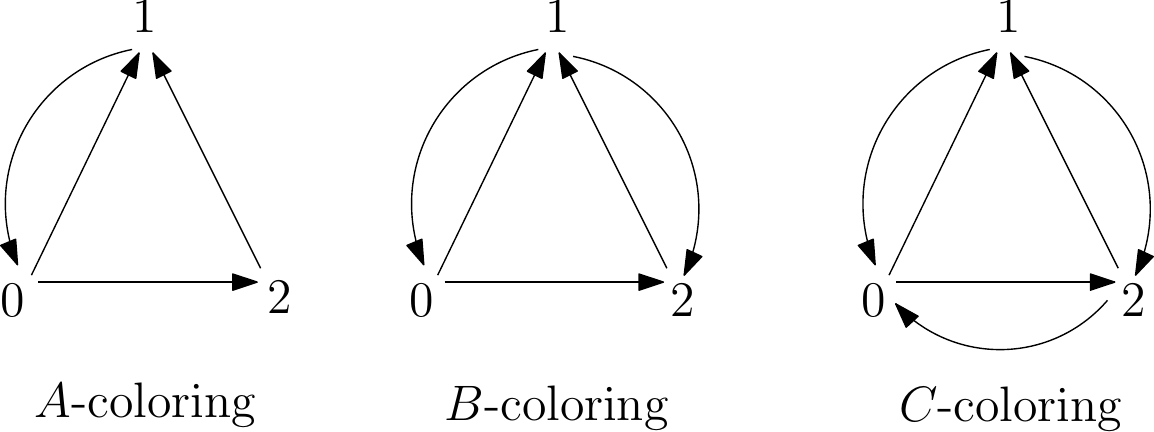}
\caption{The digraphs $A, B,$ and $C$}
\label{fig1}
\end{figure}

In contrast, for digraphs $H$, the boundary between easy and hard $H$-coloring problems is not well-understood. Some partial results have been published in \cite{MaurerSudboroughWelzl,JensenHellMac,HellZhu,Kozik} and in several other papers, but it is still not known whether all directed $H$-coloring problems are polynomial or NP-complete. This statement is in fact equivalent to the well-known Dichotomy Conjecture of CSP \cite{fv}. Even less is known about the complexity of $H$-coloring for digraphs with bounded degrees, and its relation to the conjecture of Feder, Hell, and Huang mentioned above. Of course, $K_3$ can be viewed as a symmetric digraph, namely the digraph $C$ in Figure \ref{fig1}, and in this context $C$-coloring of a digraph $G$ is precisely a $3$-coloring of the underlying undirected graph of $G$. The digraphs $A, B,$ and $C$, in Figure \ref{fig1}, are the smallest three digraphs $H$ with NP-complete $H$-coloring problems. All other digraphs $H$ with three vertices have polynomial-time solvable $H$-coloring problems \cite{JensenHellMac,MaurerSudboroughWelzl}. (All digraphs $H$ with fewer than three vertices have $H$-coloring problems that can be solved in polynomial time by $2$-SAT \cite{hombook}.) We investigate these first three interesting cases. We confirm that the $H$-coloring problems for these three digraphs, $H=A, B,$ and $C$, are polynomial time solvable for low degree bounds, and become NP-complete again when the degree bounds are high enough.

Interestingly, we can determine very precisely the highest allowable degrees that ensures the polynomiality of these problems. Let $\Delta^{+}$ denote the maximum out-degree and $\Delta^{-}$ the maximum in-degree in a digraph. We show that digraphs with $\Delta^{+} \leq 1$, $\Delta^{-} \leq 2$, or with $\Delta^{+} \leq 2$, $\Delta^{-} \leq 1$, have polynomial time $A$-coloring, $B$-coloring, and $C$-coloring algorithms, but for digraphs with $\Delta^{+} \leq 2$, $\Delta^{-} \leq 2$, all three problems are again NP-complete. The NP-completeness of $C$-coloring will imply the theorem of Holyer \cite{Holyer} mentioned earlier.

For the polynomial cases, we provide algorithms for the more general problem of {\em list $H$-coloring}. 

\section{General Digraphs}

\begin{figure}[ht]
\centering
\subfigure[Variable gadget $U$]{
\includegraphics[scale = 0.45]{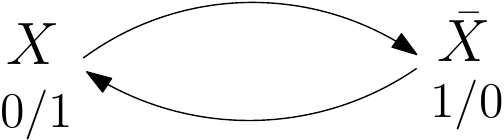}
\label{fig2a}
}
\hspace{0.3 in}
\subfigure[Clause gadget $W$]{
\includegraphics[scale = 0.45]{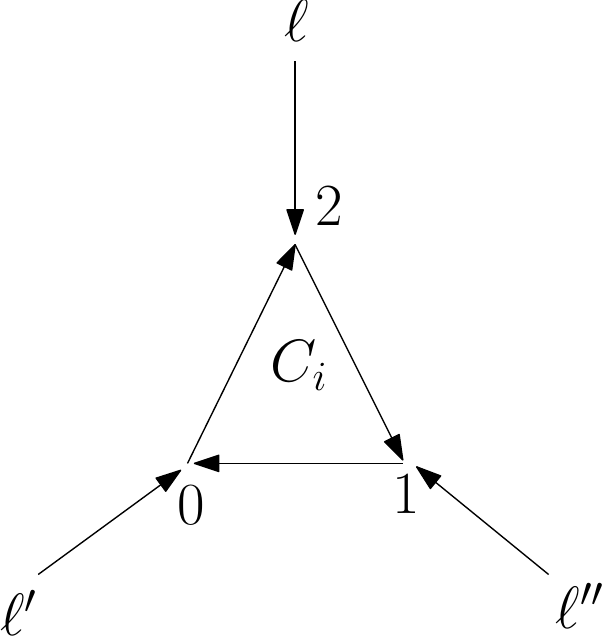}
\label{fig2b}
}
\hspace{0.4 in}
\subfigure[Combining $U$ and $W$: clause $C_i = (X_1 \vee X_2 \vee \overline{X_3})$]{
\includegraphics[scale = 0.45]{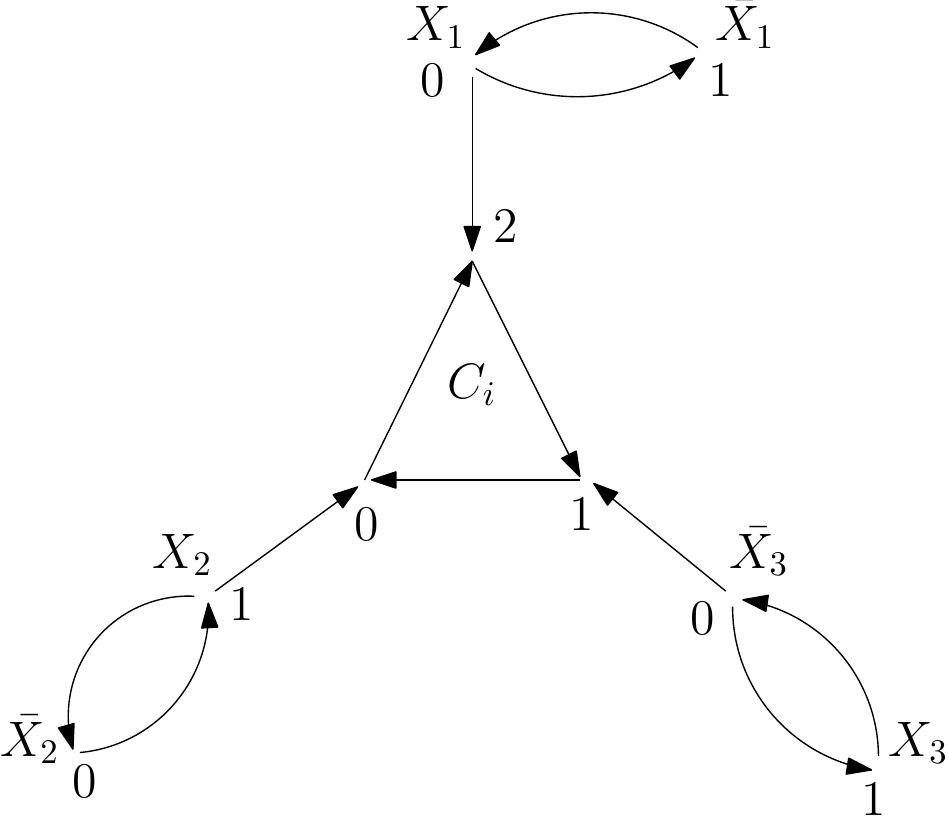}
\label{fig2c}
}
\label{fig2}
\caption[A-col]{NP-completeness of $A$-coloring}
\end{figure}

In this section, we will present new NP-completeness proofs for the $A$-coloring, $B$-coloring, and $C$-coloring of general digraphs, i.e. without any degree restrictions. It is a folklore result that these three problems are NP-complete \cite{JensenHellMac,MaurerSudboroughWelzl}. The proofs we provide here will facilitate extensions to graphs with bounded degrees.

As is standard in such proofs, we will have a {\em variable gadget} for each variable and a {\em clause gadget} for each clause. We will reduce the problem 1-in-3-SAT to the $A$- and $B$-coloring problems, i.e., for each formula $\phi$ produce a graph $G_{\phi}$ in such a way that $\phi$ is satisfiable if and only if $G_{\phi}$ has an $H$-coloring. For the $C$-coloring problem, we will show a reduction from the 3-SAT problem.

We begin with $A$-coloring.
We construct a reduction from 1-in-3-SAT as follows. Let $\phi$ be a given 3-CNF formula. Consider the variable gadget $U$ shown in Figure \ref{fig2a}. For each variable $X$ in $\phi$, we have one such gadget, with the two endpoints of $U$ corresponding to the variable $X$ and its negation $\overline{X}$. 

Figure \ref{fig2b} depicts the clause gadget $W$, for a generic clause $C_i = (\ell \vee \ell' \vee \ell'')$ with three literals $\ell, \ell', \ell''$. The vertices $\ell,\ell',\ell''$ are identified with the same literals of the corresponding variable gadgets, producing a graph $G_{\phi}$. Note that an endpoint of a variable gadget has an outgoing edge for each occurence of the corresponding variable (or its negation) in a clause of $\phi$. 

The graph $A$ has only one digon, with the vertices $0$ and $1$. Thus, there are exactly two ways of $A$-coloring the variable gadget $U$, with $X$ colored either $0$ or $1$ and $\overline{X}$ colored $1$ or $0$, respectively, as depicted in Figure \ref{fig2a}. The inner $3$-cycle of the clause gadet $W$ requires three different colors, and in fact the colors must appear in the cyclic order 0, 2, 1, as depicted in the figure, up to the symmetry of $W$. Consider now the possible colors of the vertices $\ell,\ell',\ell''$ (up to symmetry):

\begin{itemize}
 \item (0,0,0) : as $00$ is not an edge, this would not allow a $0$ on the inner cycle;
 \item (1,1,1) : as $11$ is not an edge, this would not allow a $1$ on the inner cycle;
 \item (0,1,1) : as $10$ is the only edge from $1$, this would force two $0$'s on the inner cycle;
 \item (0,0,1) : Figure \ref{fig2c} depicts the unique $A$-coloring
\end{itemize}

\begin{figure}[ht]
\centering
\subfigure[Variable gadget $V$]{
\includegraphics[scale = 0.55]{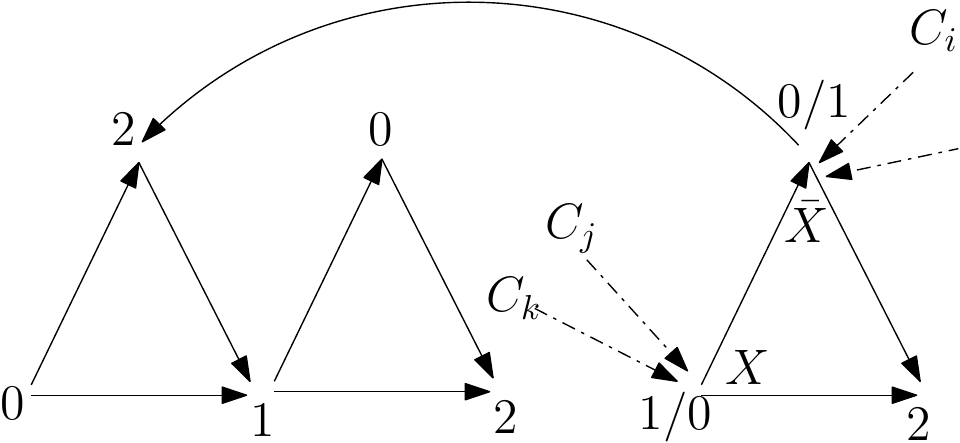}
\label{fig3a}
}
\hspace{0.7 in}
\subfigure[Clause gadget $\hat{W}$]{
\includegraphics[scale = 0.55]{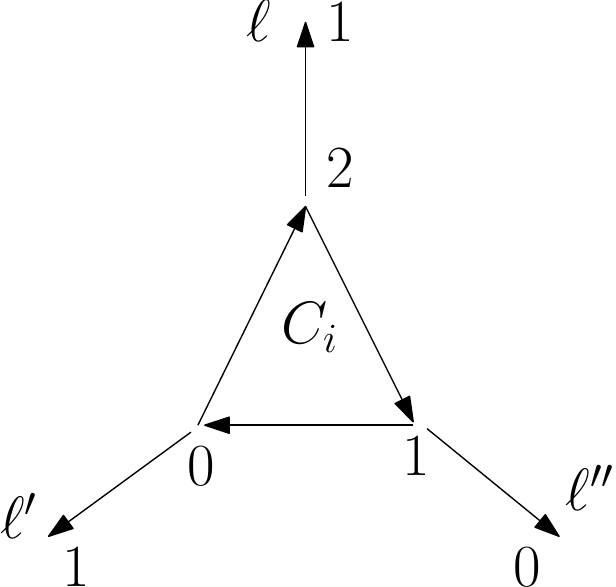}
\label{fig3b}
}
\label{fig3}
\caption[B-col]{NP-completeness of $B$-coloring}
\end{figure}

It follows that in any $A$-coloring of $G_{\phi}$, exactly one literal in each clause is colored by $1$. Thus if we assign the value True to the literals colored $1$ and value False to the literals colored $0$, we obtain a satisfying truth assignment for $\phi$. Conversely, if we start with a truth assignment and color all vertices of the variable gadgets $0$ if the literal is False and $1$ if the literal is True, there is, in each clause, exactly one literal that is colored $1$, and so the above analysis shows that the colors can be extended to the inner cycles of all the clause gadgets. In conclusion, $\phi$ is satisfiable if and only if $G_{\phi}$ is $A$-colorable.

Since the $A$-coloring problem is clearly in NP, we have the following theorem.

\begin{thm}
The $A$-coloring problem is NP-complete.
\end{thm}

A similar analysis applies for $B$-coloring. The variable gadget $V$ is shown in Figure \ref{fig3a}. It is more complex than for $A$-coloring, but there are still two vertices corresponding to the variable $X$ and its negation $\overline{X}$. It is important to observe that $V$ can be $B$-colored so that $X, \overline{X}$ are colored by $0, 1$ or by $1, 0$, but no other pair of colors. To see this, note that the vertex lying in the two adjoining triangles must be colored $1$ in any $B$-coloring of $V$, thus forcing the colors of the adjoining triangles of $V$ as depicted in the figure. Therefore, the color of $\overline{X}$ is either $0$ or $1$ (since it has an out-neighbor colored $2$). On the other hand, the color of $X$ cannot be $2$, since $X$ has two adjacent out-neighbors. Thus coloring $\overline{X}$ by $0$ forces $X$ to be colored $1$ and vice versa, and no other pair of colors for $X, \overline{X}$ is possible.

Figure \ref{fig3b} depicts the clause gadget $\hat{W}$ for $C_i = (\ell \vee \ell' \vee \ell'')$ of the 3-CNF formula. This gadget is similar to the gadget $W$ used for $A$ coloring, except the edges now point away from the center cycle, into the variable gadgets. In any $B$-coloring, the inner cycle requires three different colors. Again, colors in the order of 0, 2, 1, as depicted, give the only $B$-coloring for such a cycle (up to symmetry). The vertices $\ell,\ell',\ell''$ are identified with the corresponding literals. They can only be colored $0$ or $1$. This again gives rise to four possibilities for the colors of $\ell,\ell',\ell''$ (up to symmetry), but here only (0,1,1) is a valid choice:

\begin{enumerate}
 \item (0,0,0) : as $00$ is not an edge, this would not allow a $0$ on the inner cycle;
 \item (1,1,1) : as $11$ is not an edge, this would not allow a $1$ on the inner cycle;
 \item (0,0,1) : as $20$ is not an edge, this would force two $1$'s on the inner cycle;
 \item (0,1,1) : Figure \ref{fig3b} depicts the unique $B$-coloring.
\end{enumerate}
Now, if we assign a truth value of True to the literals that are colored $0$, and False to the literals that are colored $1$, we have a satisfying assignment for the 1-in-3-SAT problem, and vice-versa. Thus, we have the following theorem.

\begin{thm}
The $B$-coloring problem is NP-complete.
\end{thm}

Finally, we handle the case of $C$-coloring.

\begin{figure}[ht]
\centering
\subfigure[Variable gadget $T$]{
\includegraphics[scale = 0.45]{./A-col1.pdf}
\label{fig4a}
}
\hspace{0.7 in}
\subfigure[Clause gadget ${W'}$]{
\includegraphics[scale = 0.45]{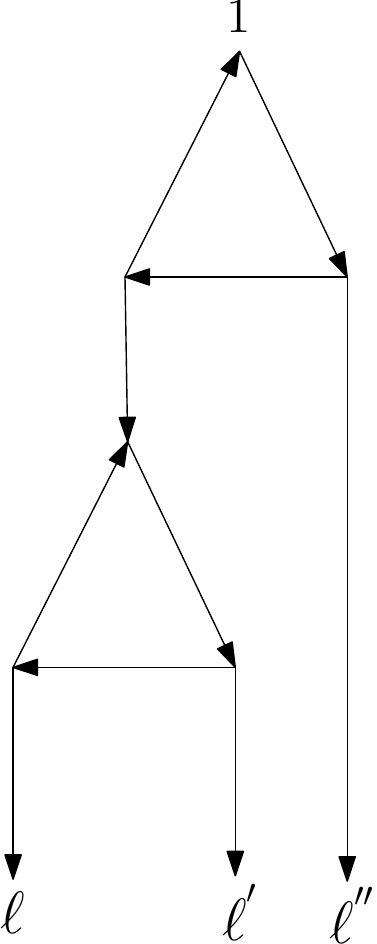}
\label{fig4b}
}
\label{fig4}
\caption[C-col]{NP-completeness of $B$-coloring}
\end{figure}

Recall that a $C$-coloring of a digraph $G$ is precisely a $3$-coloring of the underlying undirected graph of $G$. Thus we can use the proof from \cite{gj}, which constructs a reduction from 3-SAT using the clause gadget depicted in Figure \ref{fig4b}. The variable gadget depicted in Figure \ref{fig4a}, and the identification of the literals $\ell, \ell', \ell''$ with the same literals of the appropriate clause gadgets of $G_{\phi}$ are made as before. Assuming the top vertex of the clause gadget is colored $1$ as depicted, the three vertices labelled $\ell, \ell', \ell''$ (corresponding to the literals) cannot be colored by $0, 0, 0$, but can be colored by any other combination of $0, 1$. Now it is easy to see that an instance $\phi$ of $3$-SAT is satisfiable if and only if the corresponding digraph $G_{\phi}$ (i.e., its underlying undirected graph) is $3$-colorable cf. \cite{gj}.

\begin{thm}
The $C$-coloring problem is NP-complete.
\end{thm}

\section{Bounded Degree Digraphs}

We observed in the previous section that the $A$-, $B$-, and $C$-coloring problems are NP-complete for general digraphs. Now, we investigate the same problems restricted to the class of digraphs with bounded in- and out-degrees.

We begin by presenting a polynomial-time $H$-coloring algorithm for digraphs $H$ with degrees $\Delta^{+} \leq 2$, $\Delta^{-} \leq 1$, or $\Delta^{+} \leq 1$, $\Delta^{-} \leq 2$.

\begin{figure}[ht]
\centering
\subfigure[Binary tree]{
\includegraphics[scale = 0.5]{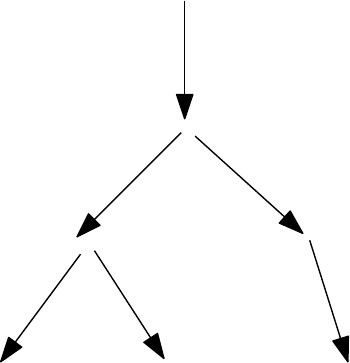}
\label{fig5a}
}
\hspace{0.5 in}
\subfigure[Directed cycle]{
\includegraphics[scale = 0.5]{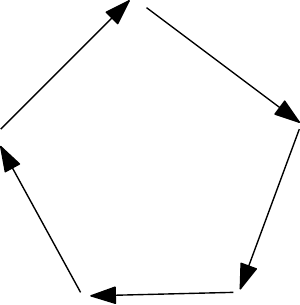}
\label{fig5b}
}
\hspace{0.5 in}
\subfigure[Cycle with hanging trees]{
\includegraphics[scale = 0.45]{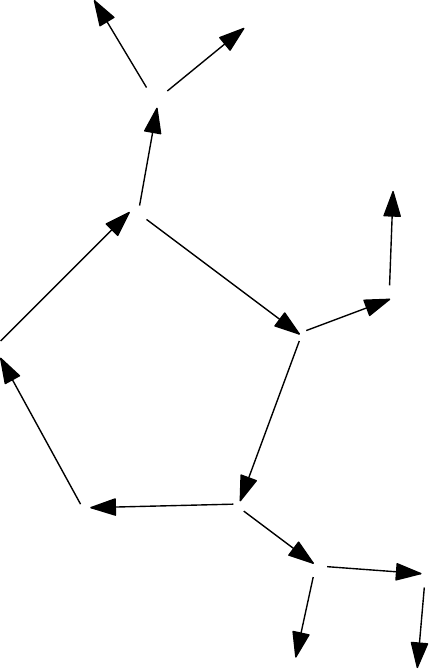}
\label{fig5c}
}
\label{fig5}
\caption[Poly]{Bounded digraphs: $\Delta^{+} \leq 2$, $\Delta^{-} \leq 1$}
\end{figure}

Let $H$ be any fixed digraph. A digraph $G$ with degrees $\Delta^{+} \leq 2$, $\Delta^{-} \leq 1$, or $\Delta^{+} \leq 1$, $\Delta^{-} \leq 2$ is either a binary tree or a directed cycle with binary trees attached at vertices of the cycle (see Figure \ref{fig5c} illustrating the case $\Delta^{+} \leq 2$, $\Delta^{-} \leq 1$). Oriented cycles and multiple directed cycles in a single component are not possible since they would violate the degree requirement at some vertex.

As mentioned earlier, we actually solve a more general problem called list $H$-coloring. To begin with, let us introduce it formally. Given a fixed digraph $H$, the problem of {\em list $H$-coloring} is the following:

\begin{center}
\textit{Given an input graph $G$, and for each $v \in V(G)$, a list $L(v) \subseteq V(H)$, is there a\\ homomorphism $f : G \to H$ such that $f(v) \in L(v)$ for all vertices $v \in V(G)$ $?$}
\end{center}

We note that by setting all lists $L(v)=V(H)$, we reduce the $H$-coloring problem to the list $H$-coloring problem. Thus, solving the list $H$-coloring problem in polynomial time also solves the $H$-coloring problem in polynomial time.

\paragraph{Algorithm.}
{\em Arc Consistency} is a basic technique from artificial intelligence, useful for solving list homomorphism problems \cite{hombook}. In particular, it consists of considering an edge $uv$ of the input graph $G$ and reducing the lists $L(u), L(v)$ so that for each $x \in L(u)$ there exists a $y \in L(v)$ with $xy \in E(H)$, and for each $z \in L(v)$ there exists a $w \in L(u)$ with $wz \in E(H)$. This kind of constraint propagation will make the entire graph arc consistent by repeating this removal as long as possible, while the lists in $G$ are changing. Since $H$ is fixed, after linearly many updates we obtain final lists satisfying all these constraints. It is easy to observe that if there is at least one empty final list, there is no homomorphism from $G$ to $H$ satisfying the original lists. 

When $G$ is a tree and the final lists are all non-empty, there is a list homomorphism $G \to H$ satisfying the original lists - simply choose one element from the final list of one fixed vertex and propagate this choice to all other vertices using the fact that both constraints were satisfied. This is not true, however, in general. For instance, when $H$ is the $2$-cycle with edges $01, 10$ and $G$ is a directed cycle of odd length (see Figure \ref{fig5b}) with all lists equal to $\{0,1\}$, then the conditions are satisfied but there is no homomorphism. If we choose one vertex of the $5$-cycle and map it to $0$, say, then the constraints will propagate around the cycle and arrive requiring the choice $1$ at the chosen vertex (and vice versa). However, if at least one vertex of a cycle has a list of size one, the constraint propagation around the cycle will work properly.

\vspace{3mm}

With this knowledge in mind, we propose our list $H$-coloring algorithm:
\begin{enumerate}
 \item Run an arc-consistency algorithm over all the edges of $G$.
 \item If some vertex has an empty list after the arc-consistency process, then, there is no list homomorphism.
 \item Else, find an $H$-coloring for each weak component of $G$ separately.
 \subitem If the component is a tree, choose one vertex and map it to one member of its list and propagate this choice to the entire component.
 \subitem If the component is a directed cycle, choose one vertex $u$, and consider $|L(u)|$ subproblems, in which the list $L(u)$ is reduced to a single vertex of $H$, i.e. for each of the possible $L(u)$ choices. Perform constraint propagation through the component based on this choice. If at least one choice leads to a homomorphism, then we have a solution, otherwise there is no solution.
 \subitem If the component has a directed cycle with hanging trees, follow the procedure to find a homomorphism for the directed cycle first. Then, extend it to find a homomorphism for each of the hanging trees.
\end{enumerate}

\begin{lem}
Let $H$ be a fixed digraph. The list $H$-coloring problem is polynomial time solvable for digraphs with $\Delta^{+} \leq 2$, $\Delta^{-} \leq 1$, or $\Delta^{+} \leq 1$, $\Delta^{-} \leq 2$.
\end{lem}

In particular, this applies to $H=A, B, C$.

\vspace{3mm}

\begin{figure}[ht]
\centering
\subfigure[Bounded degree variable gadget $U'$]{
\includegraphics[scale=0.45]{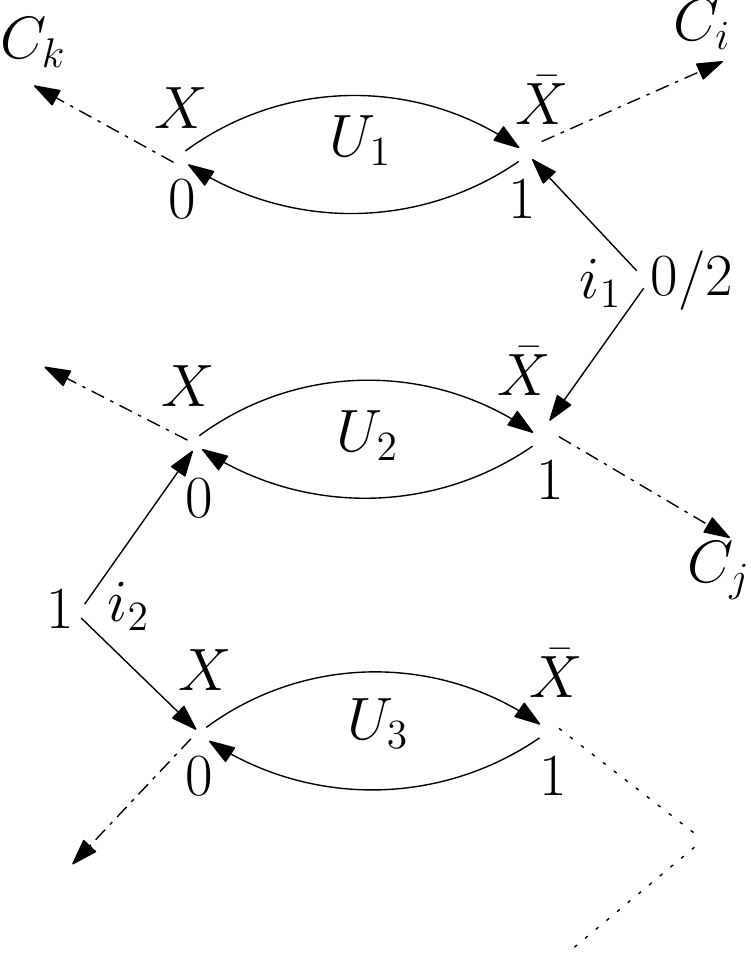}
\label{fig6a}
}
\hspace{0.7 in}
\subfigure[Bounded degree variable gadget $V'$]{
\includegraphics[scale=0.45]{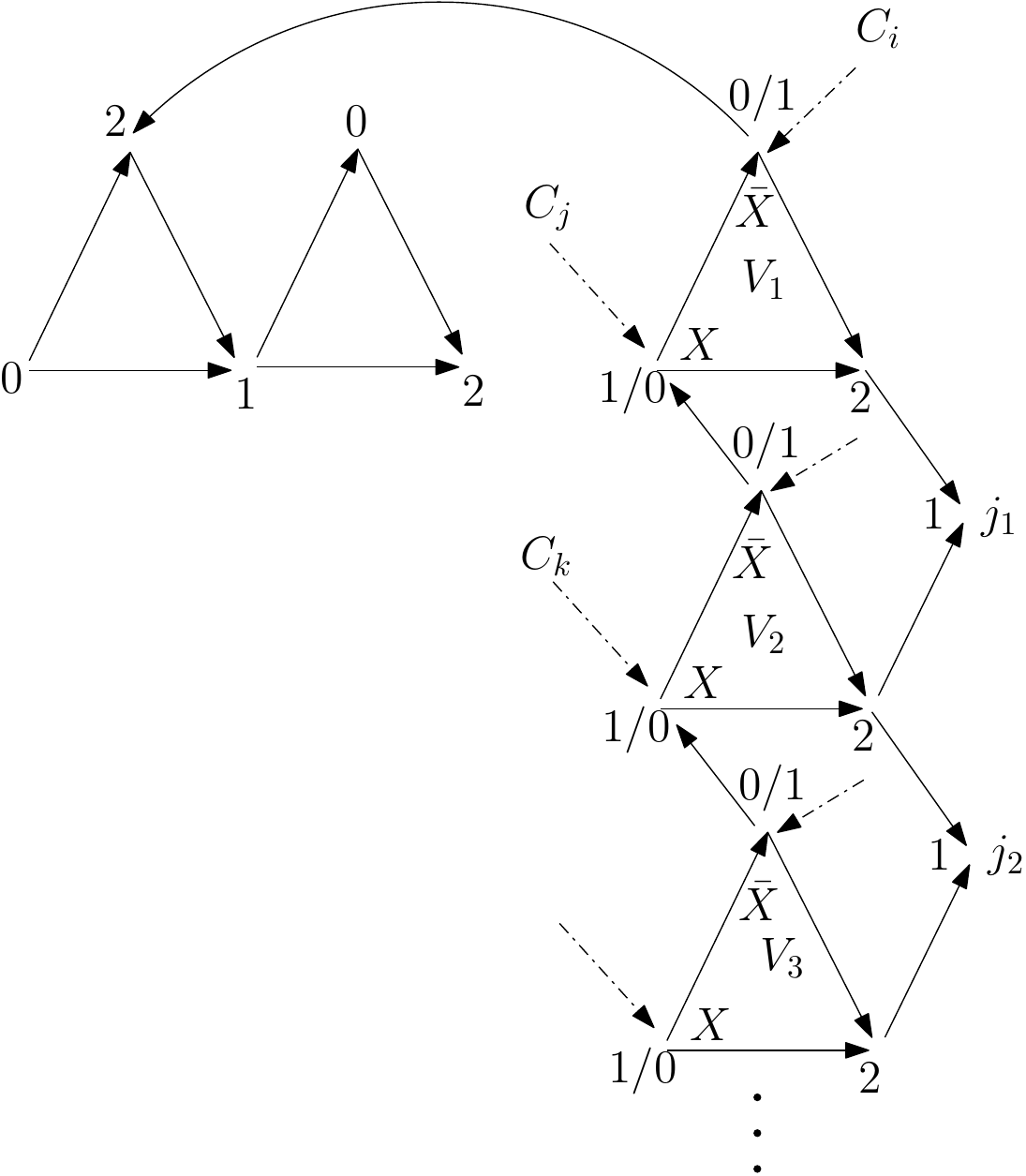}
\label{fig6b}
}
\label{fig6}
\caption[NPCnA]{$A$- and $B$-coloring for $\Delta^{+} \leq 2, \Delta^{-} \leq 2$ $\in$ NP-complete}
\end{figure}

We now focus on input digraphs with $\Delta^{+} \leq 2$, $\Delta^{-} \leq 2$. In order to prove the NP-completeness of this bounded degree setting, we need to make modifications to the gadgets used in the unbounded degree case, as the variable gadgets used in Section 2 could potentially have very high outdegree or indegree at the literal vertices. So, to keep the degree constraints satisfied, we need to construct multiple copies of each literal to join the clause gadgets of various clauses in which the literal occurs. At the same time, we have to ensure that the color assigned to the literal vertex remains same in all copies, so that we have a consistent truth value assignment.

\begin{lem}
The $A$-coloring problem is NP-complete, even when restricted to digraphs with $\Delta^{+} \leq 2$, $\Delta^{-} \leq 2$.
\end{lem}

Consider the variable gadget $U$ used in Figure \ref{fig2a}. If a literal occurs $x$ times in $\phi$, then the NP-completeness construction gives the corresponding vertex in $U$ an out-degree of $x+1$. Therefore, we shall consider a modified variable gadget $U'$, see Figure \ref{fig6a}, which has multiple copies of each literal vertex. If $X$ occurs $x$ times in $\phi$, and $\overline{X}$ occurs $x'$ times in $\phi$, then we construct $U'$ using digons $U_1, U_2, \dots, U_k$, where $k=\max(x,x')$.

As before, in any $A$-coloring each digon must be colored by $0$ and $1$. Because of the vertex $i_1$, the first two copies of $\overline{X}$ must be colored by the same color, $0$ or $1$. (Indeed, there is no vertex in $A$ that has an edge to both $0$ and $1$.) Similarly, the vertex $i_2$ ensures that the second and third copy of $X$ are colored by the same color. Repeating this argument we conclude that all vertices corresponding to $X$ have the same color, and similarly for $\overline{X}$. This allows us to apply the previous proof to a new graph $G_{\phi}$ in which each literal vertex in $U'$ is adjacent to at most one clause gadget $W$ (the same clause gadget as before). Thus $G_{\phi}$ has all in- and out-degrees bounded by $2$. It follows from the above remarks that each $A$-coloring of $G_{\phi}$ yields a satisfying truth assignment of $\phi$, and it is also easy to see that any satisfying truth assignment of $\phi$ can be extended to an $A$-coloring of $G_{\phi}$.

The analysis is again similar for $B$-coloring.

\begin{figure}[ht]
\centering
\subfigure[3-coloring of $K_{1,1,3}$]{
\includegraphics[scale=0.4]{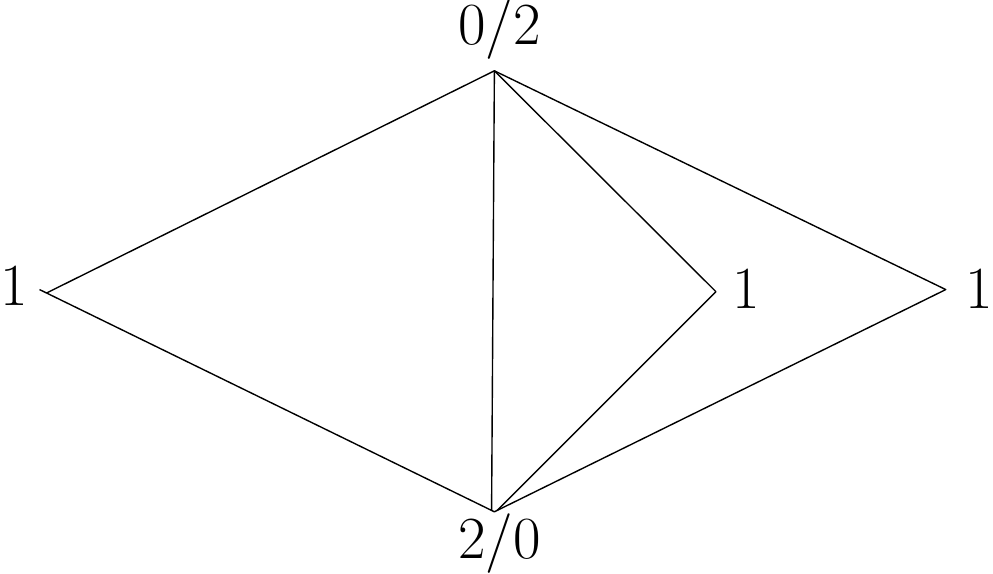}
\label{fig7a}
}
\hspace{0.7 in}
\subfigure[Bounded degree variable gadget $T'$]{
\includegraphics[scale=0.45]{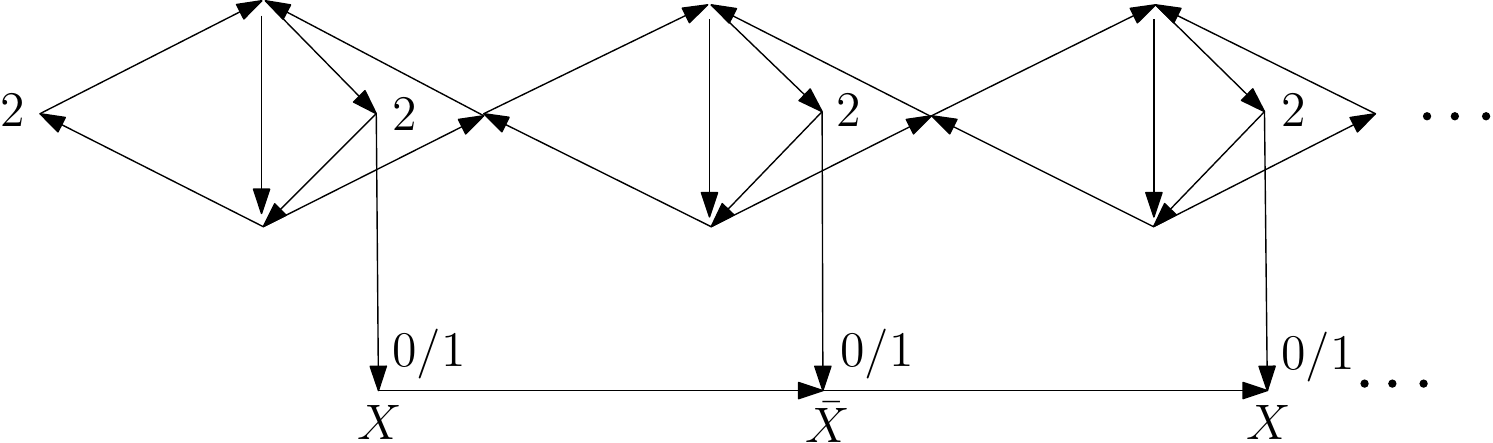}
\label{fig7b}
}
\label{fig7}
\caption[Cbound]{$C$-coloring for $\Delta^{+} \leq 2, \Delta^{-} \leq 2$ $\in$ NP-complete}
\end{figure}

\begin{lem}
The $B$-coloring problem is NP-complete, even when restricted to digraphs with $\Delta^{+} \leq 2$, $\Delta^{-} \leq 2$.
\end{lem}

In this case, we consider construct $G_{\phi}$ from the modified variable gadget $V'$ in Figure \ref{fig6b}, and the same clause gadget $\hat{W}$ as before. If $X$ occurs $x$ times in $\phi$, and $\overline{X}$ occurs $x'$ times in $\phi$, then we construct $V'$ using triangles $V_1, V_2, \dots, V_k$, where $k=\max(x,x')$, allowing for a sufficient number of copies of $X$ and $\overline{X}$ to keep the in- and out-degrees bounded by $2$. It remains to verify that the repeated copies of $X$ and $\overline{X}$ must obtain the same color in any $B$-coloring of $G_{\phi}$. As noted earlier, the third vertex of $V_1$ must be colored $2$, and hence the vertex $j_1$ must be colored $1$. This means that the third vertex of $V_2$ is colored $0$ or $2$; but $0$ is not possible as the vertex has in-degree two and $20$ is not an edge of $B$. Thus the third vertex of each $V_i$ is $2$. It now follows that in any $B$-coloring of $V'$ all copies of $X$ (and all copies of $\overline{X}$) receive the same color, $0$ or $1$. Therefore, we can repeat the previous proof of NP-completeness.

\vspace{3mm}
Finally, we discuss $C$-coloring.

\begin{lem}\label{hol}
The $C$-coloring problem is NP-complete, even when restricted to digraphs with $\Delta^{+} \leq 2$, $\Delta^{-} \leq 2$.
\end{lem}

\begin{figure}
\centering
\includegraphics[scale = 0.4]{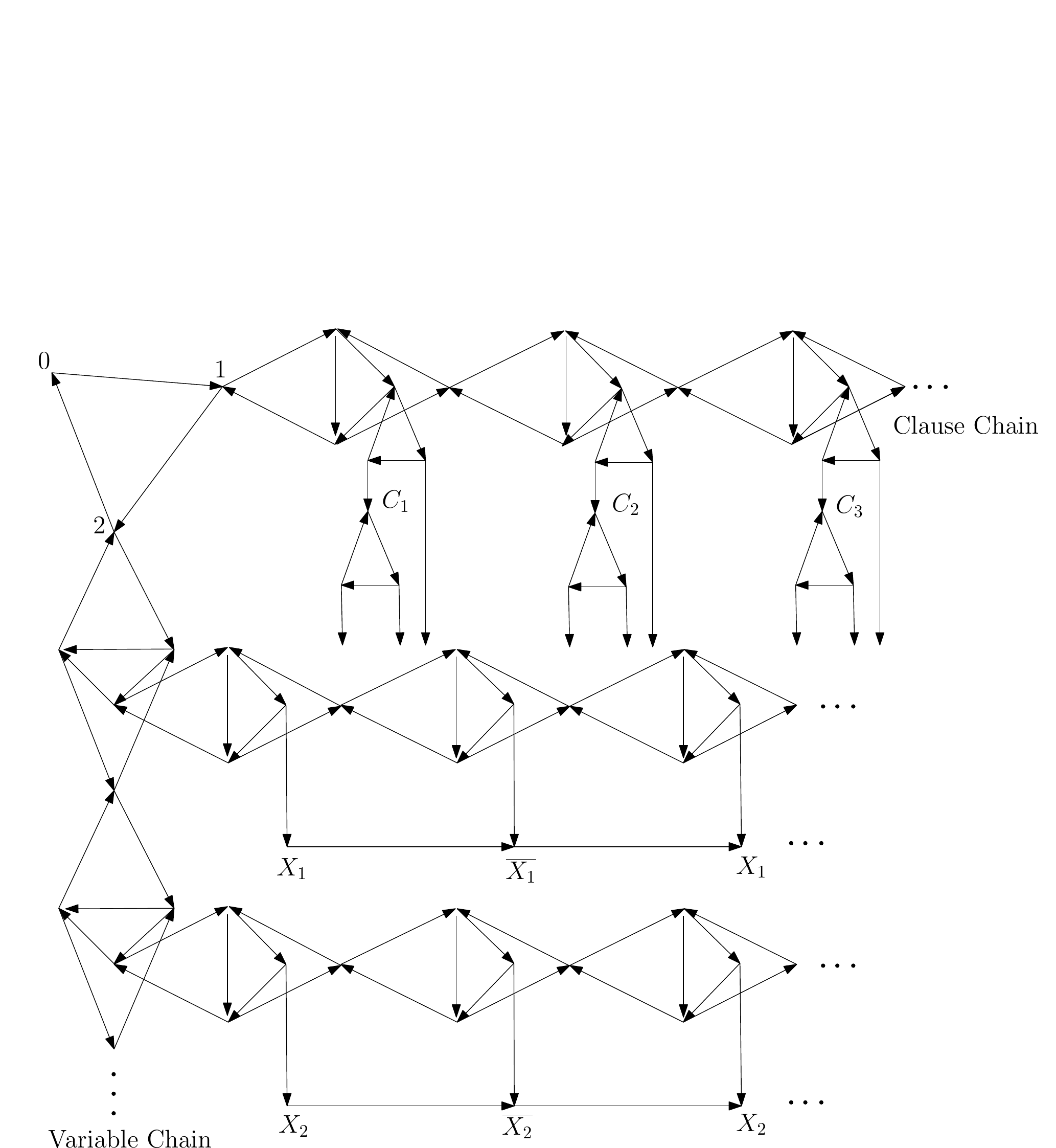}
\caption{The overall NP-completeness construction of $G_{\phi}$}
\label{fig8}
\end{figure}

Here, we use the variable gadget depicted in Figure \ref{fig7b}. It is an orientation of a chain of copies of $K_{1,1,3}$, and $K_{1,1,3}$ is a graph that takes the same color for every vertex of degree two. Thus assuming the leftmost vertex is colored $2$ as depicted, all the vertices corresponding to $X$ and $\overline{X}$ are colored by $0$ and $1$ or $1$ and $0$ respectively. The overall construction is depicted in Figure \ref{fig8}. There is a basis triangle from which the construction emanates, and which is assumed to be colored by $0, 1, 2$ as depicted (by renaming the colors if necessary). There is a {\em clause chain} (of copies of $K_{1,1,3}$) providing us with a sufficient number of vertices colored $1$ for each clause $C_i$ of $\phi$. We use the same clause gadget $W'$ as in the unbounded setting. Finally, there is a {\em variable chain} (of copies of $K_{1,1,3}$) yielding a sufficient number of vertices colored $2$ for each variable $X_j$. The literal vertices of the clause gadgets are identified with the corresponding literals in the variable chains as before. It is easy to check that the maximum in- and out-degree is two, and that $\phi$ is satisfiable if and only if $G_{\phi}$ is $3$-colorable.

Combining all the results for the restricted setting of bounded degree input graphs, we have our main theorem.

\begin{thm}
Deciding $A$-, $B$-, and $C$-colorability is in $\mathbb{P}$ for input digraphs with $\Delta^{+} \leq 2, \Delta^{-} \leq 1$ or $\Delta^{+} \leq 1, \Delta^{-} \leq 2$, and is NP-complete for digraphs with $\Delta^{+} \leq 2, \Delta^{-} \leq 2$.
\end{thm}

\section{Conclusions}

There are three smallest digraphs $H$ with NP-complete $H$-coloring problems, namely $A, B,$ and $C$ from Figure \ref{fig1}. We have shown that for all three, the $H$-coloring problem remains NP-complete if the in- and out-degrees are bounded from above by $2$, but become polynomial-time solvable if either of these degree bounds is lowered to $1$. This suggests that the conjecture of Feder, Hell, and Huang, claiming that NP-complete coloring problems remain NP-complete even for bounded degree graphs provided the degree bound is high enough, may also hold for digraphs. This is the first result of this type. Moreover, we were successful at delineating very precisely the boundary between the degree bounds that imply polynomiality and those that do not.

We conclude by noting that our proof implies the following classical result. The proof is exactly the same as for Lemma \ref{hol}, except the edges of the gadgets are not oriented.

\begin{thm}\cite{Holyer}
The $3$-coloring problem is NP-complete even when restricted to graphs with maximum degree $4$.
\end{thm}

\section*{Acknowledgements}
The authors thank the MITACS Globalink Internship Program for undergraduate students for funding the internship of the second author which made this collaboration possible. The authors would also like to acknowledge the hospitality of Simon Fraser Univeristy, and especially of the IRMACS center, where most of this was research was done, during the research internship of the second author.

\bibliographystyle{plain}
\bibliography{H_coloring_aurosish_pavol}
\end{document}